\documentclass[10pt]{article}
\input{xypic}
\usepackage{amsmath}
\usepackage{amssymb}
\usepackage{latexsym}
\usepackage{enumerate}

\newtheorem{theorem}{Theorem}[section]
\newtheorem{proposition}[theorem]{Proposition}
\newtheorem{corollary}[theorem]{Corollary}
\newtheorem{lemma}[theorem]{Lemma}

\newtheorem{remark}[theorem]{Remark}

\newenvironment{proof}{\noindent{\sc Proof}}{\hfill\qed}
\def\endrem{}

\newcommand{\F}{{\mathbb F}}

\newcommand{\dst}{\displaystyle}

\newcommand{\Hom}{{\rm Hom}}
\newcommand{\Aut}{{\rm Aut}}
\newcommand{\Out}{{\rm Out}}

\newcommand{\id}{{\rm Id}}
\newcommand{\qed}{\quad\lower0.05cm\hbox{$\Box$}}
\newcommand{\AG}{{\Aut(G)}}
\newcommand{\AH}{{\Aut(H)}}
\newcommand{\OG}{{\Out(G)}}
\newcommand{\OH}{{\Out(H)}}

\newcommand{\ls}[2]{{\,^{#1}\!#2}}                    

\newcommand{\downarrowright}[1]{\downarrow
\rlap{\raise0.1cm\hbox{$\scriptstyle{#1}$}}}

\newcommand{\downarrowleft}[1]{\rlap{\kern-0.2cm
\raise0.1cm\hbox{$\scriptstyle{#1}$}}\downarrow}

\newcommand{\uparrowright}[1]{\uparrow
\rlap{\lower0.1cm\hbox{$\scriptstyle{#1}$}}}

\newcommand{\uparrowleft}[1]{\rlap{\kern-0.2cm
\lower0.1cm\hbox{$\scriptstyle{#1}$}}\uparrow}

\newcommand{\ra}{\rightarrow}
\newcommand{\lra}{\longrightarrow}

\newcommand{\mono}{\hookrightarrow}
\newcommand{\lemono}{\hookleftarrow}

\def\rmono{\rto|<\hole|<<\ahook}

\def\umono{\ar@{_{(}->}[u]}
\def\uumono{\ar@{_{(}->}[uu]}
\def\dmono{\dto|<\hole|<<\ahook}
\def\ddmono{\ddto|<\hole|<<\ahook}
\def\lmono{\ar@{_{(}->}[l]}
\def\llmono{\ar@{_{(}->}[ll]}

\begin{document}
\title{Finite simple groups and localization}

\author{
{\sc Jos\'e L. Rodr\'{\i}guez
\thanks{Partially supported by DGESIC grant PB97-0202 and
the Swiss National Science Foundation.}} \\[2pt]
{\em Departamento de Geometr\'\i a, Topolog\'\i a y Qu\'\i mica Org\'anica}\\
{\em Universidad de Almer\'\i a, E--04120 Almer\'\i a, Spain}\\
{\em e-mail: {\tt jlrodri@ual.es}} \\[1cm]
{\sc J\'er\^ome Scherer\thanks{Partially supported by the Swiss National Science Foundation.}
 and Jacques Th\'evenaz}\\[2pt]
{\em Institut de Math\'ematiques, Universit\'e de Lausanne}\\
{\em CH--1015 Lausanne, Switzerland}\\
{\em e-mail: {\tt jerome.scherer@ima.unil.ch}}\\
{\em e-mail: {\tt jacques.thevenaz@ima.unil.ch}}
}

\date{}
\maketitle

\begin{abstract}
The purpose of this paper is to explore the concept of localization,
which comes from homotopy theory, in the context of finite simple groups.
We give an easy criterion for a finite simple group to be a localization of
some simple subgroup and we apply it in various cases.
Iterating this process allows us to connect many simple groups by a sequence of
localizations.
We prove that all sporadic simple groups (except possibly the Monster) and several
groups of Lie type are connected to alternating groups.
The question remains open whether or not there are several connected components
within the family of finite simple groups.
\end{abstract}


\section*{Introduction}

The concept of localization plays an important role in homotopy theory.
The introduction by Bousfield of homotopical localization functors in \cite{Bou76}
and more recently its popularization by Farjoun in \cite{Farj96} has
led to the study of localization functors in other categories.
Special attention has been set on the category of groups $Gr$, as the effect
of a homotopical localization on the fundamental group is often best
described by a localization functor ${\rm L}\colon Gr \to Gr$.

A {\bf localization functor\/} is a pair $({\rm L}, \eta)$
consisting of a functor ${\rm L}\colon Gr \to Gr$ together with a natural
transformation $\eta\colon \id \to {\rm L}$ from the identity functor, such that $L$
is idempotent, meaning that the two morphisms
$\eta_{{\rm L}G}, {\rm L}(\eta_G)\colon {\rm L}G \to {\rm L}{\rm L}G$ coincide
and are isomorphisms.
A group homomorphism $\varphi\colon H\to G$ is called in turn a {\bf localization\/}
if there exists a localization functor $({\rm L}, \eta)$ such that
$G = {\rm L}H$ and $\varphi=\eta_H \colon H \to {\rm L}H$ (but we note that the
functor ${\rm L}$ is not uniquely determined by~$\varphi$).
In this situation, we often say that $G$ is a localization of~$H$.
A very simple characterization of localizations can be given without mentioning
localization functors:
A group homomorphism $\varphi\colon H\to G$ is a localization if and only if $\varphi$
induces a bijection
\addtocounter{theorem}{1}
\begin{equation}
\label{bijection}
\varphi^*\colon \Hom(G,G)\cong \Hom(H, G)
\end{equation}
as mentionned in \cite[Lemma~2.1]{Cas00}.
In the last decade several authors (Casacuberta, Farjoun, Libman, Rodr\'\i guez) have
directed their efforts towards deciding which algebraic properties are
preserved under localization.
An exhaustive survey about this problem is nicely exposed in \cite{Cas00} by Casacuberta.
For example, any localization of an abelian group is again abelian.
Similarly, nilpotent groups of class at most~2 are preserved, but the question remains open
for arbitrary nilpotent groups.
Finiteness is not preserved, as shown by the example $A_n\to SO(n-1)$ for $n \geq 10$
(this is the main result in \cite{Lib98}).
In fact, it has been shown by  G\"obel and Shelah in~\cite{GS01} that any non-abelian
finite simple group has arbitrarily large localizations (a previous version of this result,
assuming the generalized continuum hypothesis, was obtained in \cite{GRS99}).
In particular it is not easy to determine all possible localizations
of a given object.
Thus we restrict ourselves to the study of finite groups and
wonder if it would be possible to understand the finite localizations of a given
finite simple group.
This paper is a first step in this direction.

Libman \cite{Lib99} observed recently that the inclusion
$A_n\mono A_{n+1}$ of alternating groups is a localization if~$n\geq 7$.
His motivation was to find a localization where new torsion elements appear
(e.g. $A_{10}\mono A_{11}$ is such a localization since $A_{11}$ contains elements
of order 11).
In these examples, the groups are simple, which simplifies considerably
the verification of formula (\ref{bijection}).
It suffices to check if $\Aut(G)\cong \Hom(H,G) - \{ 0\}.$

\medskip

This paper is devoted to the study of the behaviour of injective localizations
with respect to simplicity.
We first give a criterion for an inclusion of a simple group in a finite simple
group to be a localization.
We then find several infinite families of such localizations, for example $L_2(p)
\mono A_{p+1}$ for any prime $p \geq 13$ (cf. Proposition~\ref{linear2}).
Here $L_2(p)=PSL_2(p)$ is the projective special linear group.
It is striking to notice that the three conditions that appear in our criterion for
an inclusion of simple groups $H \mono G$ to be a localization already appeared in the
literature.
For example the main theorem of~\cite{KW90a} states exactly that
$J_3 \mono E_6(4)$ is a localization (see Section~\ref{section rigid}).
Similarly the main theorem in \cite{SWW00} states that $Sz(32) \mono E_8(5)$ is
a localization.
Hence the language of localization theory can be useful to shortly reformulate some rather
technical properties.

By Libman's result, the alternating groups $A_n$, for $n\geq 7$, are all
connected by a sequence of localizations.
We show that $A_5\mono A_6$ is also a localization.
A more curious way allows us to connect $A_6$ to $A_{7}$ by a zigzag of localizations:
$$
A_6 \mono T \mono Ru \lemono L_2(13) \mono A_{14} \lemono \cdots
\lemono A_7
$$
where $T$ is the Tits group, and $Ru$ the Rudvalis group.
This yields the concept of rigid component of a simple group.
The idea is that among all inclusions $H \mono G$, those that are localizations
deserve our attention because of the ``rigidity condition" imposed by~(\ref{bijection}):
Any automorphism of $G$ is completely determined by its restriction to~$H$.
So, we say that two groups $H$ and $G$ lie in the same {\bf rigid component\/} if $H$
and $G$ can be connected by a zigzag of inclusions which are all localizations.

Many finite simple groups can be connected to the alternating groups.
Here is our main result:

\bigskip\noindent
{\bf Theorem} {\it The following finite simple groups all lie in the same rigid
component:
\renewcommand{\baselinestretch}{.8}
\small\normalsize
\begin{itemize}
\item[(i)] All alternating groups $A_n$ ($n\geq 5$).
\item[(ii)] The Chevalley groups $L_2(q)$ where $q$ is a prime power~$\geq 5$.
\item[(iii)] The Chevalley groups $U_3(q)$ where $q$ is any prime power.
\item[(iv)] The Chevalley groups $G_2(p)$ where $p$ is an odd prime such that $(p+1, 3) = 1$.
\item[(v)] All sporadic simple groups, except possibly the Monster.
\renewcommand{\baselinestretch}{1.2}
\small\normalsize
\item[(vi)] The Chevalley groups $L_3(3)$, $L_3(5)$, $L_3(11)$, $L_4(3)$,
$U_4(2)$, $U_4(3)$, $U_5(2)$, $U_6(2)$, 
$S_4(4)$, $S_6(2)$, $S_8(2)$, $D_4(2)$, $\ls{2}{D_4(2)}$, $\ls{2}{D_5(2)}$, $\ls{3}{D_4(2)}$,
$D_4(3)$, $G_2(2)'$, $G_2(4)$, $G_2(5)$, $G_2(11)$,
$E_6(4)$, $F_4(2)$, and $T=\ls{2}{F_4(2)'}$.
\end{itemize}
}

The proof is an application of the localization criteria which are given in
Sections~\ref{section criterion} and \ref{section alternating}, but
requires a careful checking in the ATLAS \cite{atlas}, or in the more
complete papers about maximal subgroups of finite simple groups (e.g.
\cite{Kleidman88}, \cite{LW91}, \cite{WilsonBaby}).
We do not know if the Monster can be connected to the alternating groups.

It is still an open problem to know how many rigid components of finite simple groups
there are, even though our main theorem seems to suggest that there is only one.
We note that the similar question for non-injective localizations has a trivial
answer (see Section~1).

Let us finally mention that simplicity is not necessarily preserved by localization.
This is the subject of the separate paper \cite{RSV}, where we exhibit for example
a localization map from the Mathieu group $M_{11}$ to the double cover of
the Mathieu group~$M_{12}$.
In our context this implies that the rigid component of a simple group may contain a
non-simple group.
This answers negatively a question posed both by Libman in \cite{Lib99} and Casacuberta in
\cite{Cas00} about the preservation of simplicity.
In these papers it was also asked whether perfectness is preserved, but we leave this
question unsolved.

\bigskip

{\it Acknowledgments:}
We would like to thank Antonio Viruel, Jean Michel, as well
as the referee for helpful comments.
The first-named author also thanks the Institut de Math\'ematiques
de l'Universit\'e de Lausanne for its kind hospitality and fantastic
views of the Alps.

\section{A localization criterion}
\label{section criterion}

Let us fix from now on a finite simple group~$G$.
In Theorem~\ref{maincriterion} below we list necessary and sufficient
conditions for an inclusion $H \mono G$
between two non-abelian finite simple groups to be a
localization.
These conditions are easier to deal with
if the groups $H$ and $G$ satisfy some extra assumptions, as
we show in the corollaries after the theorem.
The proof is a variation of that of Corollary~4 in \cite{GRS99}.

We note here that we only deal with injective group homomorphisms because
non-injective localizations abound.
For example, for any two finite groups $G_1$ and $G_2$ of coprime orders,
$G_1\times G_2\to G_1$ and  $G_1\times G_2\to G_2$ are localizations.
So the analogous concept of rigid component defined using
non-injective localizations has no interest, since obviously any two finite
groups are in the same component.

If the inclusion $i\colon H \mono G$ is a localization, then so is the inclusion
$H' \mono G$ for any subgroup $H'$ of $G$ which is isomorphic to~$H$.
This shows that the choice of the subgroup $H$ among isomorphic
subgroups does not matter.

Let  $c\colon G \to \AG$ be the natural injection of $G$ defined as
$c(g)=c_g\colon G\to G$, where $c_g$ is the inner automorphism
given by $x\mapsto gxg^{-1}$.
We shall always identify in this way a simple group $G$ with a subgroup of~$\AG$
and the quotient $\OG = \AG/G$ is then called the {\bf group of outer automorphisms}
of~$G$.

\begin{lemma}
\label{conjugateautomorphism}
Let $G$ be a non-abelian simple group.
Then the following diagram commutes
\[
\diagram
G \rto^\alpha \dto_c & G \dto^c \cr
\AG \rto_{c_\alpha} & \AG
\enddiagram
\]
for any automorphism $\alpha \in \AG$.
\end{lemma}

\begin{proof}
This is a trivial check.
\end{proof}

\begin{lemma}
\label{extendautomorphism}
Let $H$ be a non-abelian simple subgroup of a finite simple group~$G$.
Suppose that the inclusion $i\colon H \mono G$ extends to an inclusion of their
automorphism groups $i\colon \AH \mono \AG$,
i.e., the following diagram commutes
\[
\diagram
H \rmono^i \dto_c & G \dto^c \cr
\AH \rmono_i & \AG
\enddiagram
\]
Then every automorphism $\alpha\colon H \ra H$ extends to an automorphism
$i(\alpha)\colon G \ra G$.
\end{lemma}

\begin{proof}
We have to show that the following square commutes:
\[
\diagram
H \rto^\alpha \dmono_i & H \dmono^i \cr
G \rto_{i(\alpha)} & G
\enddiagram
\]
To do so we consider this square as the left-hand face of the cubical diagram
\[
\diagram
H \rrto^c \ddmono \drto_\alpha && \AH
\ddto|<\hole|<<\ahook|(0.5)\hole \drto^{c_\alpha} \cr
& H \ddmono \rrto^c && \AH \ddmono \cr
G \drto_{i(\alpha)} \rrto|(0.42)\hole^c && \AG \drto^{c_{i(\alpha)}} \cr
& G \rrto_c && \AG
\enddiagram
\]
The top and bottom squares commute by Lemma~\ref{conjugateautomorphism}.
The front and back squares are the same and commute by assumption.
The right-handsquare commutes as well because $i$ is a homomorphism.
This forces the left-hand square to commute and we are done.
\end{proof}

\begin{remark}
\label{converse}
{\rm
As shown by the preceding lemma, it is stronger to require that the inclusion
$i\colon {H}\mono{G}$ extends to an inclusion $i\colon \AH \mono \AG$ than to require that
every automorphism of $H$ extends to an automorphism of~$G$.
In general we have an exact sequence
$$
1 \ra C_\AG(H) \lra N_\AG(H) \lra \AH
$$
so the second condition is equivalent to the fact that this is a short exact sequence.
However, in the presence of the condition $C_\AG(H) = 1$, which
plays a central role in this paper, we find that $N_\AG(H) \cong \AH$.
Thus any automorphism of $H$ extends to a unique automorphism of $G$, and this
defines a homomorphism $i\colon \AH \mono \AG$ extending the inclusion $H \mono G$.
Therefore, if the condition $C_\AG (H) = 1$ holds, we have a converse of the above
lemma and both conditions are equivalent.
We will use the first in the statements of the following results, even though it is the
stronger one.
It is indeed easier to check in the applications.
}
\end{remark}
\endrem

\begin{theorem}
\label{maincriterion}
Let $H$ be a non-abelian simple subgroup of a finite simple group~$G$ and let
$i\colon H \mono G$ be the inclusion.
Then $i$ is a localization if and only if the following three conditions are satisfied:
\begin{enumerate}
\item The inclusion $i\colon H \mono G$ extends to an inclusion $i\colon \AH \mono \AG$.
\item Any subgroup of $G$ which is isomorphic to $H$ is conjugate to
$H$ in $\Aut(G)$.
\item The centralizer $C_{\Aut(G)} (H) = 1$.
\end{enumerate}
\end{theorem}

\begin{proof}
If $i$ is a localization, all three conditions have to be satisfied.
Indeed given an automorphism $\alpha \in \AH$, formula~(\ref{bijection}) tells us
that there exists a unique group homomorphism $\beta\colon G \ra G$ such that
$\beta \circ i = i \circ \alpha$.
Since $G$ is finite and simple, $\beta$ is an automorphism. 
We set then $i(\alpha) = \beta$ and condition~(1) follows.
Given a subgroup $j\colon H' \mono G$ and an isomorphism $\phi\colon H \ra H'$, a similar argument
with formula~(\ref{bijection}) ensures the existence of an automorphism $\beta \in \AG$
such that $\beta \circ i = j \circ \phi$.
Thus condition~(2) holds.
Finally condition~(3) is also valid since the unique extension of $1_H$ to $G$ is the
identity.

Assume now that all three conditions hold.
For any given homomorphism ${\varphi\colon H\to G}$, we need a unique
homomorphism $\Phi\colon G\to G$ such that $\Phi \circ i =\varphi$.
If $\varphi$ is trivial, we choose of course  the trivial homomorphism $\Phi\colon G\to G$.
It is unique since $H$ is in the kernel of $\Phi$, which must be equal to $G$ by simplicity.
Hence, we can suppose that $\varphi$ is not trivial.
Since $H$ is simple we have that $\varphi(H)\leq G$ and $H\cong \varphi(H)$.

By (2) there is an automorphism $\alpha\in \AG$
such that $c_\alpha(\varphi(H))= H$, or equivalently by
Lemma~\ref{conjugateautomorphism}, $\alpha(\varphi(H)) = H$.
Therefore the composite map
$$
H\stackrel{\varphi}{\longrightarrow}\varphi(H)
\stackrel{\alpha|_{\varphi(H)}}{\longrightarrow}H
$$
is some automorphism $\beta$ of~$H$.
By condition (1) this automorphism of $H$ extends to
an automorphism $i(\beta)\colon G \ra G$.
That is, the following square commutes:
\[
\diagram
H \rto^\beta \dmono_i & H \dmono^i \cr
G \rto_{i(\beta)} & G
\enddiagram
\]
The homomorphism $\Phi = \alpha^{-1}i(\beta)$ extends $\varphi$ as desired.
We prove now it is unique.
Suppose that $\Phi'\colon G\to G$ is a homomorphism such that $\Phi' \circ
i =\varphi$.
Then, since $G$ is simple, $\Phi' \in \AG$.
The composite $\Phi^{-1} \Phi'$ is an element in the centralizer
$C_{\Aut(G)} (H)$, which is trivial by (3).
This finishes the proof of the theorem.
\end{proof}

\begin{remark}
\label{orbit}
{\rm The terminology used for condition~(2) is that two subgroups $H$ and $H'$ in $G$
{\bf fuse} in $\AG$ if there is an automorphism $\alpha \in \AG$ such that $\alpha(H) = H'$.
Assuming condition (1) in Theorem~\ref{maincriterion} we have a short exact
sequence
$$
1 \ra C_\AG(H) \lra N_\AG(H) \lra \AH \ra 1
$$
(compare with Remark~\ref{converse}).
Moreover there are $\dst\frac{\mid\AG \mid}{\mid N_\AG(H) \mid}$ subgroups in the
conjugacy class of $H$ in $\AG$ and $\dst\frac{\mid G \mid}{\mid N_G(H) \mid}$ subgroups
in the conjugacy class of $H$ in~$G$.
If condition (2) also holds, this implies that the total number of conjugacy classes
of subgroups isomorphic to $H$ in $G$ is equal to
$$
\frac{\mid\AG \mid}{\mid N_\AG(H) \mid} \cdot \frac{\mid N_G(H) \mid}{\mid G \mid}
 = \frac{\mid\OG \mid \cdot \mid G \mid}{\mid \AH \mid \cdot \mid C_\AG(H) \mid} \cdot
\frac{\mid N_G(H) \mid}{\mid G \mid}
$$

\medskip
\noindent
Condition (3) is thus equivalent to the following one, which is sometimes easier to verify:
\begin{itemize}
\item[{\it 3'.}]
{\it The number of conjugacy classes of subgroups of $G$ isomorphic to~$H$ is equal to
$$
\dst\frac{\mid\OG \mid}{\mid\OH\mid} \cdot \frac{\mid N_G(H) \mid}{\mid H\mid}.
$$
}
\end{itemize}
}
\end{remark}
\endrem

We obtain immediately the following corollaries. 
Using the terminology in \cite[p.158]{Rot95}, recall that a group is {\bf complete}
if it has trivial centre and every automorphism is inner.

\begin{corollary}
\label{complete}
Let $H$ be a non-abelian simple subgroup of a finite simple group~$G$ and let
$i\colon H \mono G$ be the inclusion.
Assume that $H$ and $G$ are complete groups.
Then $i$ is a localization if and only if the following two
conditions are satisfied:
\begin{enumerate}
\item Any subgroup of $G$ which is isomorphic to $H$ is conjugate to~$H$.
\item $C_G(H) = 1$. \hfill{\qed}
\end{enumerate}
\end{corollary}

The condition $C_G(H) = 1$ is here equivalent to $N_G(H) = H$.
This is often easier to check.
It is in particular always the case when $H$ is a maximal subgroup of~$G$.
This leads us to the next corollary.

\begin{corollary}
\label{max}
Let $H$ be a non-abelian simple subgroup of a finite simple group~$G$ and let
$i\colon H \mono G$ be the inclusion.
Assume that $H$ is a maximal subgroup of~$G$.
Then $i$ is a localization if and only if the following three
conditions are satisfied:
\begin{enumerate}
\item The inclusion $i\colon H \mono G$ extends to an inclusion $i\colon \AH \mono \AG$.
\item Any subgroup of $G$ which is isomorphic to $H$ is conjugate to
$H$ in~$\Aut(G)$.
\item There are $\dst{\frac{|\Out(G)|}{|\Out(H)|}}$ conjugacy classes of subgroups
isomorphic to $H$ in~$G$.
\end{enumerate}
\end{corollary}

\begin{proof}
Since $H$ is a maximal subgroup of $G$, $N_G(H) = H$.
The corollary is now a direct consequence of Theorem~\ref{maincriterion} taking into
account Remark~\ref{orbit} about the number of conjugacy classes of subgroups
of $G$ isomorphic to~$H$.
\end{proof}

\section{Localization in alternating groups}
\label{section alternating}

We describe in this section a method for finding localizations of finite simple
groups in alternating groups.
Let $H$ be a simple group and $K$ a subgroup of index~$n$.
The (left) action of $H$ on the cosets of $K$ in $H$ defines a
{\bf permutation representation\/} $H \to S_n$ as in~\cite[Theorem~3.14, p.53]{Rot95}.
The {\bf degree\/} of the representation is the number $n$ of cosets.
As $H$ is simple, this homomorphism is actually an inclusion $H \mono A_n$.
Recall that $\Aut(A_n)=S_n$ if $n\geq 7$.

\begin{theorem}
\label{largestmaximal}
Let $H$ be a non-abelian finite simple group and $K$ a maximal subgroup
of index~$n\geq 7$.
Suppose that the following two conditions hold:
\begin{enumerate}
\item
The order of $K$ is maximal (among all maximal subgroups).
\item
Any subgroup of $H$ of index $n$ is conjugate to~$K$.
\end{enumerate}
Then the permutation representation $H \mono A_n$ is a localization.
\end{theorem}

\begin{proof}
We show that the conditions of Theorem~\ref{maincriterion} are satisfied,
starting with condition~(1).
Since $K$ is maximal, it is self-normalizing and therefore the action of $H$
on the cosets of $K$ is isomorphic to the conjugation action of $H$ on the set of
conjugates of~$K$.
By our second assumption, this set is left invariant under~$\AH$.
Thus the action of~$H$ extends to~$\AH$ and this yields the desired
extension $\AH\to S_n=\Aut(A_n)$.

To check condition (2) of Theorem~\ref{maincriterion}, let $H'$ be a
subgroup of $A_n$ which is isomorphic to~$H$ and denote by $\alpha\colon H\to H'$
an isomorphism.
Let $J$ be the stabilizer of a point in $\{1,\ldots,n\}$
under the action of~$H'$.
Since the orbit of this point has cardinality~$\leq n$, the index of~$J$ is
at most~$n$, hence equal to~$n$ by our first assumption.
Thus $H'$ acts transitively.
So $H$ has a second transitive action via~$\alpha$ and the action of~$H'$.
For this action, the stabilizer of a point is a subgroup of $H$ of index~$n$,
hence conjugate to~$K$ by assumption.
So $K$ is also the stabilizer of a point for this second action and this shows
that this action of $H$ is isomorphic to the permutation action of $H$ on the
cosets of~$K$, that is, to the first action.
It follows that the permutation representation
$H \stackrel{\alpha}{\longrightarrow} H' \mono A_n$ is conjugate in
$S_n$ to~$H \mono A_n$.

Finally, since $H$ is a transitive subgroup of $S_n$ with maximal
stabilizer, the centralizer $C_{S_n}(H)$ is trivial by \cite[Theorem~4.2A (vi)]{DiMo}
and thus condition~(3) of Theorem~\ref{maincriterion} is satisfied.
\end{proof}

\bigskip

Among the twenty-six sporadic simple groups, twenty have a subgroup which
satisfies the conditions of Theorem~\ref{largestmaximal}.

\begin{corollary}
\label{twentysporadic}
The following inclusions are localizations:

\noindent
${M_{11} \mono A_{11}\,}$,
${M_{22} \mono A_{22}\,}$,
${M_{23} \mono A_{23}\,}$,
${M_{24} \mono A_{24}\,}$,
${J_1 \mono A_{266}\,}$,
${J_2 \mono A_{100}\,}$,

\noindent
${J_3 \mono A_{6156}\,}$,
${J_4 \mono A_{173067389}\,}$,
${HS \mono A_{100}\,}$,
${McL \mono A_{275}\,}$,
${Co_1 \mono A_{98280}\,}$,

\noindent
${Co_2 \mono A_{2300}\,}$,
${Co_3 \mono A_{276}\,}$,
${Suz \mono A_{1782}\,}$,
${He \mono A_{2058}\,}$,
${Ru \mono A_{4060}\,}$,

\noindent
${Fi_{22} \mono A_{3510}\,}$,
${Fi_{23} \mono A_{31671}\,}$,
${HN \mono A_{1140000}\,}$,
${Ly \mono A_{8835156}\,}$.
\end{corollary}

\begin{proof}
In each case, it suffices to check in the ATLAS~\cite{atlas} that the
conditions of Theorem~\ref{largestmaximal} are satisfied.
It is however necessary to check the complete list of maximal subgroups in~\cite{KPW89}
for the Fischer group $Fi_{23}$ and~\cite{KW88} for the Janko group~$J_4$.
\end{proof}

\bigskip

We obtain now two infinite families of localizations.
The classical projective special linear groups $L_2(q) = PSL_2(q)$ of type $A_1(q)$,
as well as the projective special unitary groups $U_3(q) = PSU_3(q)$ of type $\ls{2}{A_2(q)}$,
are almost all connected to an alternating group by a localization.
Recall that the notation $L_2(q)$ is used only for the simple projective special linear groups,
that is if the prime power~$q \geq 4$.
Similarly the notation $U_3(q)$ is used for~$q > 2$.

\bigskip

\begin{proposition}
\label{linear2}
$(i)$ The permutation representation $L_2(q) \mono A_{q+1}$ induced by the action of
$SL_2(q)$ on the projective line is a localization for any prime power $q \not\in
\{4, 5, 7, 9, 11 \}$.

\hskip 2.5cm $(ii)$ The permutation representation $U_3(q) \mono A_{q^3+1}$ induced
by the action of $SU_3(q)$ on the set of isotropic points in the projective plane is a
localization for any prime power~$q \neq 5$.
\end{proposition}

\begin{proof}
We prove both statements at the same time.
The group $L_2(q)$ acts on the projective line, whereas $U_3(q)$ acts on the set of
isotropic points in the projective plane.
In both cases, let $B$ be the stabilizer of a point for this action (Borel subgroup).
Let us also denote by $G$ either $L_2(q)$ or $U_3(q)$,
where $q$ is a prime power as specified above, and $r$ is $q+1$, or $q^3 + 1$ respectively.
Then $B$ is a subgroup of $G$ of index $r$ by \cite[Satz~II-8.2]{Hu} and
\cite[Satz~II-10.12]{Hu}.

By \cite[Satz~II-8.28]{Hu}, which is an old theorem of Galois
when $q$ is a prime, the group $L_2(q)$ has no non-trivial permutation representation
of degree less than~$r$ if $q \not\in \{4, 5, 7, 9, 11 \}$.
The same holds for $U_3(q)$ by~\cite[Table~1]{cooperstein} if~$q \neq 5$.
Thus $B$ satisfies condition (1) of Theorem~\ref{largestmaximal}.

It remains to show that condition (2) is also satisfied.
The subgroup $B$ is the normalizer of a Sylow $p$-subgroup $U$, and $B = UT$,
where $T$ is a complement of $U$ in~$B$.
If $N$ denotes the normalizer of $T$ in $G$, we know that $G = UNU$.
This is the Bruhat decomposition (for more details see \cite[Chapter~8]{MR53:10946}).
We are now ready to prove that any subgroup of $G$ of index $r$ is conjugate to~$B$.
Let $H$ be such a subgroup.
It contains a Sylow $p$-subgroup, and we can thus assume it actually contains~$U$.
Since $G$ is generated by $U$ and $N$, the subgroup $H$ is generated by $U$ and~$N \cap H$.
Assume $H$ contains an element $x \in N - T$.
The Weyl group $N/T$ is cyclic of order two, generated by the class of $x$
(the linear group $L_2(q)$ is a Chevalley group of type $A_1$ and the unitary group $U_3(q)$
is a twisted Chevalley group of type $\ls{2}{A_2}$).
Moreover $G = \; <U, xUx^{-1}>$ by \cite[Theorem~2.3.8~(e)]{98j:20011}
and both $U$ and its conjugate $xUx^{-1}$ are contained in~$H$.
This is impossible because $H \neq G$, so $N \cap H = T \cap H$.
It follows that $H$ is contained in $<U, T> \;  = B$.
But $H$ and $B$ have the same order and therefore~$H = B$.
\end{proof}

\begin{remark}
{\rm This proof does not work for the action of $L_{n+1}(q)$ on the
$n$-dimensional projective space if $n\geq 2$, because there is a second action
of the same degree, namely the action on the set of all hyperplanes in
$(\F_q)^{n+1}$.
Thus there is another conjugacy class of subgroups of the same
index, so condition (1) does not hold.}
\end{remark}
\endrem

\section{Proof of the main theorem}
\label{section rigid}

In order to prove our main theorem, we have to check that any group of the
list is connected to an alternating group by a zigzag of localizations.
When no specific proof is indicated for an inclusion to be a
localization, it means that all the necessary information for checking
conditions (1)-(3) of Theorem~\ref{maincriterion} is available in the
ATLAS~\cite{atlas}.
By $C_2$ we denote the cyclic group of order~2.


\bigskip
\noindent
{\it (i) Alternating groups}.


The inclusions $A_n \mono A_{n+1}$, for $n\geq 7$, studied by Libman in
\cite[Example~3.4]{Lib99} are localizations by Corollary~\ref{max},
with $\Out(A_n) \cong C_2 \cong \Out(A_{n+1})$.
The inclusion $A_5 \mono A_6$ is a localization as well, since
we have $\Out(A_6) \cong (C_2)^2$, $\Out(A_5) \cong C_2$, and
there are indeed two conjugacy classes of subgroups of $A_6$
isomorphic to $A_5$ with fusion in~$\Aut(A_6)$.
The inclusion $A_6 \mono A_7$ is not a localization,
but we can connect these two groups via a zigzag of localizations,
for example as follows:
$$
A_6 \mono T \mono Ru \lemono L_2(13) \mono A_{14}
$$
where $T$ denotes the Tits group, $Ru$ the Rudvalis group and the
last arrow is a localization by Proposition~\ref{linear2}.


\bigskip
\noindent
{\it (ii)  Chevalley groups $L_2(q)$}.


By Proposition~\ref{linear2}, all but five linear groups $L_2(q)$ are
connected to an alternating group.
The groups $L_2(4)$ and $L_2(5)$ are isomorphic to $A_5$,
and $L_2(9) \cong A_6$.
We connect $L_2(7)$ to $A_{28}$ via a chain of two localizations
$$
L_2(7) \mono U_3(3) = G_2(2)' \mono A_{28}
$$
where we use Theorem~\ref{largestmaximal} for the second map.
Similarly, we connect $L_2(11)$ to $A_{22}$ via the Mathieu group $M_{22}$, using
Corollary~\ref{twentysporadic}:
$$
L_2(11) \mono M_{22} \mono A_{22}.
$$


\bigskip
\noindent
{\it (iii)  Chevalley groups $U_3(q)$}.


For $q \neq 5$, we have seen in Proposition~\ref{linear2}~(ii) that
$U_3(q) \mono A_{q^3 + 1}$ is a localization.
One checks in \cite[p.34]{atlas} that there is a localization $A_7 \mono U_3(5)$,
which connects $U_3(5)$ to the alternating groups.


\bigbreak
\noindent
{\it (iv)  Chevalley groups $G_2(p)$}.


When $p$ is an odd prime such that $(p+1, 3) = 1$, we will see in
Proposition~\ref{unitary3} that $U_3(p) \mono G_2(p)$ is a localization.
We can conclude by (iii), since 5 is not a prime in the considered family.


\bigskip
\noindent
{\it (v)  Sporadic simple groups}.

By Corollary~\ref{twentysporadic}, we already know that twenty sporadic
simple groups are connected with some alternating group.
We now show how to connect all the other sporadic groups, except the Monster
for which we do not know what happens.

For the Mathieu group~$M_{12}$, we note that the inclusion $M_{11} \mono
M_{12}$ is a localization because there are two conjugacy classes of
subgroups of $M_{12}$ isomorphic to $M_{11}$ (of index~11) with fusion in
$\Aut(M_{12})$ (cf. \cite[p.33]{atlas}).
We conclude by Corollary~\ref{max}.

The list of all maximal subgroups of $Fi_{24}'$ is given in \cite{LW91}
and one applies Theorem~\ref{maincriterion} to show that $He \mono Fi_{24}'$ is a
localization (both groups have $C_2$ as group of outer automorphisms).

Looking at the complete list of maximal subgroups of the Baby Monster
$B$ in \cite{WilsonBaby}, we see that $Fi_{23} \mono B$ is a localization,
as well as $Th \mono B$, $HN \mono B$, and $L_2(11) \mono B$ (see
Proposition~4.1 in \cite{WilsonBaby}).
This connects Thompson's group $Th$ and the Baby Monster (as well as
the Harada-Norton group $HN$) to the Fischer groups and also to the Chevalley
groups~$L_2(q)$.

Finally we consider the O'Nan group~$O'N$.
By~\cite[Proposition~3.9]{WilsonON} we see that $M_{11} \mono O'N$ satisfies
conditions (1)-(3) of Corollary~\ref{max} and thus is a localization.


\bigbreak
\noindent
{\it (vi)  Other Chevalley groups}.

The construction of the sporadic group $Suz$ provides a
sequence of five graphs (the Suzuki chain) whose groups of automorphisms are
successively $\Aut(L_2(7))$, $\Aut(G_2(2)')$, $\Aut(J_2)$, $\Aut(G_2(4))$
and $\Aut(Suz)$ (see \cite[p.108-9]{gore}).
Each one of these five groups is an extension of $C_2$ by the appropriate
finite simple group.
All arrows in the sequence
\addtocounter{theorem}{1}
\begin{equation}
\label{Suzuki}
L_2(7) \mono G_2(2)' \mono J_2 \mono G_2(4) \mono Suz
\end{equation}
are thus localizations by Corollary~\ref{max} because they are actually
inclusions of the largest maximal subgroup (cf. \cite{atlas}).
This connects the groups $G_2(2)'$ and $G_2(4)$ to alternating groups since we already
know that $Suz$ is connected to~$A_{1782}$ by Corollary~\ref{twentysporadic}.
Alternatively, note that $G_2(4)\lemono L_2(13)\mono A_{14}$ are localizations,
using Proposition~\ref{linear2} for the second one.

The Suzuki group provides some more examples of localizations:
${L_3(3) \mono Suz}$ by \cite[Section~6.6]{WilsonSuz}, and $U_5(2) \mono Suz$
by \cite[Section~6.1]{WilsonSuz}.
We also have localizations
$$
A_9 \mono D_4(2)=O_8^+(2) \mono F_4(2) \lemono \ls{3}{D_4(2)}
$$
which connect these Chevalley groups (see Proposition~\ref{F4} for the last arrow).
We are able to connect three symplectic groups since $A_8 \mono S_6(2)$ and $S_4(4) \mono He$
are localizations, as well as $S_8(2) \mono A_{120}$ by Theorem~\ref{largestmaximal}.
This allows us in turn to connect more Chevalley groups as $U_4(2) \mono S_6(2)$, and
$O_8^{-}=\ls{2}{D_4(2)} \mono S_8(2)$ are all localizations.

Each of the following localizations involves
a linear group and connects some new group to the component of the alternating groups:

\noindent
$L_2(11) \mono U_5(2)$, $L_3(3) \mono T$, $L_2 (7) \mono L_3(11)$, and $L_4(3) \mono F_4(2)$.

The localization $U_3(3) \mono G_2(5)$ connects $G_2(5)$ and
thus $L_3(5)$ by Proposition~\ref{linear3} below.
Likewise, since we just showed above that $L_3(11)$ belongs to the same rigid component,
then so does~$G_2(11)$.

Next $M_{22} \mono U_6(2)$ and $A_{12} \mono O_{10}^- = \ls{2}{D_5(2)}$
are also localizations.

In the last three localizations, connecting the groups $U_4(3)$, $E_6(4)$, and
$D_4(3)$, the orders of the outer automorphism groups is larger than~2.
Nevertheless, Theorem~\ref{maincriterion} applies easily.
There is a localization $A_7 \mono U_4(3)$.
There are four conjugacy classes of subgroups of $U_4(3)$ isomorphic to $A_7$, all
of them being maximal.
The dihedral group $D_8 \cong \Out(U_4(3))$ acts transitively on those classes
and $S_7$ is contained in $\Aut(U_4(3))$ (see~\cite[p.52]{atlas}).

We have also a localization $J_3 \mono E_6(4)$.
Here $\Out(E_6(4)) \cong D_{12}$ and  there are exactly six conjugacy classes of sugroups
isomorphic to $J_3$ in $E_6(4)$ which are permuted transitively by~$D_{12}$.
This is exactly the statement of the main theorem of \cite{KW90a}.

Finally $D_4(2) \mono D_4(3)$ is a localization.
Here we have $\Out(D_4(2)) \cong S_3$ and  $\Out(D_4(3)) \cong S_4$.
There are four conjugacy classes of subgroups of $D_4(3)$ isomorphic to~$D_4(2)$.


\section{Other localizations}
\label{section other}


In this section, we give further examples of localizations between simple
groups.
We start with three infinite families of localizations.
Except the second family, we do not know if the groups belong to the rigid
component of alternating groups.

\begin{lemma}
\label{only}
Let $p$ be an odd prime with $(3, p-1) = 1$.
If $H \leq M \leq G_2(p)$ are subgroups with $H \cong L_3(p)$ and $M$ maximal,
then $M \cong \Aut (L_3(p))$.
\end{lemma}

\begin{proof}
The order of $L_3(p)$ is $p^3(p^3-1)(p^2-1)$, which is larger than~$p^6$.
So the main theorem in \cite{LS87} implies that any maximal sugbroup of $G_2(p)$
containing $H$ has to be one of \cite[Table~1, p.300]{LS87} or a parabolic subgroup.
Of course $| L_3(p) |$ has to divide~$| M |$.
Therefore $M$ can not be parabolic because the order of a parabolic subgroup of $G_2(p)$
is $p^6(p^2-1)(p-1)$ (see \cite[Theorem~A]{Kleidman88}).
When $p \neq 3$ there are only two maximal subgroups left to deal with.
When $p=3$, we could also have used the ATLAS \cite[p.60]{atlas}.
In both cases the only maximal subgroup whose order is divisible by $| L_3(p) |$
is isomorphic to~$\Aut (L_3(p))$.
\end{proof}

\begin{proposition}
\label{linear3}
Let $p$ be an odd prime with $(3, p-1) = 1$.
Then there is a localization $L_3(p) \mono G_2(p)$.
\end{proposition}

\begin{proof}
The strategy is to verify conditions (1)-(3) of Theorem~\ref{maincriterion}, or rather
conditions (1), (2), and (3') of Remark~\ref{orbit}.
Let us first review some facts from \cite[Proposition~2.2]{Kleidman88}.
There exist two subgroups $L_+ \leq K_+ \leq G_2(p)$ with $L_+ \cong L_3(p)$ a subgroup of
index 2 in~$K_+$.
This is the inclusion we consider here.
Moreover $K_+ = N_{G_2(p)}(L_+)$ (this is Step 3 in \cite[Proposition~2.2]{Kleidman88}), and
$N_{\Aut(G_2(p))}(L_+) \cong \Aut(L_3(p))$ (Step 1).
Thus $K_+ \cong \Aut(L_3(p))$ and condition~(1) of Theorem~\ref{maincriterion} obviously holds.
To check condition~(2) consider a subgroup $H \leq G_2(p)$ with $H \cong L_3(p)$.
By the above lemma $H$ must be contained in a maximal subgroup isomorphic to~$\Aut(L_3(p))$.
When $p \neq 3$ the group $G_2(p)$ is complete and \cite[Theorem~A]{Kleidman88} shows
that $H$ is contained in some conjugate of~$K_+$.
Therefore $H$ is conjugate to $L_+$ in~$G_2(p)$.
When $p=3$, the group $G_2(p)$ is not complete.
The group $\Out(G_2(3))$ is cyclic of order 2 as there is a ``graph automorphism"
$\gamma\colon G_2(3) \ra G_2(3)$ which is not an inner automorphism (we follow the notation from
\cite[Proposition~2.2~(v)]{Kleidman88}).
Then $K_+$ and $\gamma(K_+)$ are not conjugate in $G_2(3)$, but of course they are
in~$\Aut(G_2(3))$.
Again Theorem~A in \cite{Kleidman88} shows that $H$ is conjugate
either to $L_+$ or to $\gamma(L_+)$ in $G_2(3)$ and condition (2) is also satisfied.
Finally condition~(3') is valid since $\Out(L_3(p))$ is cyclic of order~2.
The number of conjugacy classes of subgroups isomorphic to $L_3(p)$ in $G_2(p)$ is 1
when $p \neq 3$ and 2 when~$p=3$.
\end{proof}

\begin{proposition}
\label{unitary3}
Let $p$ be an odd prime with $(3, p+1) = 1$.
Then there is a localization $U_3(p) \mono G_2(p)$.
\end{proposition}

\begin{proof}
The proof is similar to that of the preceding proposition.
Apply also \cite[Proposition~2.2]{Kleidman88} for the subgroup~$K_{-}$.
\end{proof}

\begin{proposition}
\label{F4}
There is a localization $\ls{3}{D_4(p)} \mono F_4(p)$ for any prime~$p$.
\end{proposition}

\begin{proof}
We have $\Out(F_4(2)) \cong C_2$ while for an odd prime $p$, $F_4(p)$ is complete.
On the other hand $\Out(\ls{3}{D_4(p)})$ is cyclic of order~3.
By \cite[Proposition~7.2]{LS87} there are exactly $(2, p)$ conjugacy classes of
subgroups isomorphic to $\ls{3}{D_4(p)}$ in $F_4(p)$, fused by an automorphism if~$p=2$.
Applying Theorem~\ref{maincriterion}, we see that the inclusion
$\ls{3}{D_4(p)} \mono \Aut(\ls{3}{D_4(p)}) \mono F_4(p)$ given by
\cite[Table~1, p.300]{LS87} is a localization.
\end{proof}

\bigskip

We have seen various localizations involving sporadic
groups in the proof of the main theorem.
We give here further examples.

We start with the five Mathieu groups.
Recall that the Mathieu groups $M_{12}$ and $M_{22}$ have $C_2$ as
outer automorphism groups, while the three other Mathieu groups are complete.
The inclusions $M_{11} \mono M_{23}$ and $M_{23} \mono M_{24}$ are
localizations by Corollary~\ref{complete}.
We have already seen in Section~\ref{section rigid}
that the inclusion $M_{11} \mono M_{12}$ is a localization.
The inclusion ${M_{12} \mono M_{24}}$ is also a localization.
Indeed $\Aut(M_{12})$ is the stabilizer in $M_{24}$ of a  pair of dodecads,
the stabilizer of a single dodecad is a copy of~$M_{12}$.
Up to conjugacy, these are the only subgroups of $M_{24}$ isomorphic to
$M_{12}$ and thus the condition (3') in Remark~\ref{orbit} is satisfied.
Similarly $M_{22} \mono M_{24}$ is also a localization, because
$\Aut(M_{22})$ can be identified as the stabilizer of a duad (a pair of octads) in $M_{24}$
whereas $M_{22}$ is the pointwise stabilizer (see \cite[p.39 and p.94]{atlas}).
In short we have the following diagram, where all inclusions are localizations:
\[
\diagram
M_{11} \rmono \dmono & M_{23} \dmono & \cr
M_{12} \rmono & M_{24} & M_{22} \; . \lmono
\enddiagram
\]


We consider next the sporadic groups linked to the Conway group~$Co_1$.
Inside $Co_1$ sits $Co_2$ as stabilizer of a certain vector $OA$ of type
2 and $Co_3$ as stabilizer of another vector $OB$ of type~3.
These vectors are part of a triangle $OAB$ and its stabilizer is the group $HS$,
whereas its setwise stabilizer is~$\Aut(HS)$.
The Conway groups are complete, the smaller ones are maximal simple subgroups of
$Co_1$ and there is a unique conjugacy class of each of them in $Co_1$ as indicated
in the ATLAS \cite[p.180]{atlas}.
Hence $Co_2 \mono Co_1$ and $Co_3 \mono Co_1$ are localizations by Corollary~\ref{complete}.
Likewise the inclusions $HS \mono Co_2$
and $McL \mono Co_3$ are also localizations:
They factor through their group of automorphisms, since for example $\Aut(McL)$ is the
setwise stabilizer of a triangle of type $223$ in the Leech lattice,
a vertex of which is stabilized by~$Co_3$.
Finally, $M_{22} \mono HS$ is a localization for similar reasons, since
any automorphism of $M_{22}$ can be seen as an automorphism of the
Higman-Sims graph (cf. \cite[Theorem~8.7 p.273]{beth}).
We get here the following diagram of localizations:
\[
\diagram
M_{22} \dmono & McL \rmono & Co_3 \dmono \cr
HS \rmono & Co_2 \rmono & Co_1 
\enddiagram
\]

\noindent
Some other related localizations are $M_{23} \mono Co_3$, $M_{23}
\mono Co_2$ and $M_{11} \mono HS$.


We move now to the Fischer groups and Janko's group~$J_4$.
The inclusion ${T \mono Fi_{22}}$ is a localization (both have $C_2$ as
outer automorphism groups) as well as ${M_{12} \mono Fi_{22}}$, and
${A_{10} \mono Fi_{22}}$.
Associated to the second Fischer group, we have a chain of localizations
$$
A_{10} \mono S_8(2) \mono Fi_{23}.
$$
By \cite[Theorem~1]{KPW89} the inclusion $A_{12} \mono Fi_{23}$ is also a localization.
Moreover ${M_{11} \mono J_4}$ and $M_{23} \mono J_4$ are localizations by
Corollaries~6.3.2 and 6.3.4 in~\cite{KW88}.


\bigskip

Let us list now without proofs a few inclusions we know to be
localizations.
We start with two examples of localizations of alternating groups:
$A_{12}\mono HN$, and $A_7 \mono Suz$ by
\cite[Section~4.4]{WilsonSuz}.
Finally we list a few localizations of Chevalley groups:
$L_2(8) \mono S_6(2)$, $L_2(13) \mono  G_2(3) $, $L_2(32) \mono J_4 $
(by \cite[Proposition~5.3.1]{KW88}), $U_3(3)\mono S_6(2)$,
$\ls{3}{D_4(2)} \mono Th$, $G_2(5) \mono Ly$, $E_6(2) \mono E_7(2)$, and $E_6(3) \mono E_7(3)$.
The inclusion $E_6(q) \mono E_7(q)$ is actually a localization if and only if
$q =2$ or $q=3$  by \cite[Table~1]{LS87}.
The main theorem in \cite{SWW00} states that
$Sz(32) \mono E_8(5)$ is a localization.
There is a single conjugacy class of subgroups isomorphic to $Sz(32)$ in $E_8(5)$,
and  $\Out(Sz(32))$ is cyclic of order~5.

\bigskip

We have seen a great deal of localizations of finite simple groups, and one could
think at this point that they abound in nature.
This is of course not so:
During our work on this paper, we came across many more inclusions of simple groups that
are not localizations.
They do not appear here for obvious reasons. 
On the other hand, our list is certainly far from being complete.
It would be nice for example to find other infinite families of localizations
among groups of Lie type and to determine if they are connected to the
alternating groups.
Another interesting task would be to find out which simple groups are local with
respect to a given localization $i: H \mono G$, or even better to compute the
localization of other simple groups with respect to $i$.
Will they still be finite simple groups?


\end{document}